\input amstex
\frenchspacing
\documentstyle{amsppt}
\magnification=\magstep1
\vsize=18.5cm 
\footline{\hfill\fiverm version October 28, 2013}
\def\today{\ifcase\month\or January \or February\or
March\or April\or May\or June\or July\or August\or
September\or October\or November\or December\fi
\space\number\day, \number\year}
\NoBlackBoxes
\def\FF{\bold F}

\def\Aut{\mathop{\text{\rm Aut}}\nolimits}
\def\Gal{\mathop{\text{\rm Gal}}}
\def\Hom{\mathop{\text{\rm Hom}}}

\def\isar{\ \smash{\mathop{\longrightarrow}\limits^\sim}\ }

\def\ker{\mathop{\text{\rm ker}}}
\def\im{\mathop{\text{\rm im}}}

\def\mapright#1{\ \smash{\mathop{\longrightarrow}\limits^{#1}}\ }
\def\legendre{\overwithdelims()}
\def\mod{\bmod}
\def\congr{ \equiv }
\def\iso{ \cong }
\def\tto{\longrightarrow}

\def\Que{\bold Q}

\def\CC{\bold C}
\def\Zee{\bold Z}

\def\muhat{\widehat{\mu}}
\newcount\refCount
\def\newref#1 {\advance\refCount by 1
\expandafter\edef\csname#1\endcsname{\the\refCount}}
\newref CP   
\newref CW   
\newref HO   
\newref LEa  
\newref LMS  
\newref MA   
\newref MO   
\newref PA   
\newref SC   
\newref ST   
\topmatter
\title
Computing higher rank primitive root densities
\endtitle
\author P. Moree and  P. Stevenhagen\endauthor
\address 
Mathematisch Instituut, Universiteit Leiden, Postbus 9512, 2300 RA Leiden,
The Netherlands
\endaddress
\email psh\@math.leidenuniv.nl\endemail
\address 
Max-Planck-Institut f\"ur Mathematik,
Vivatsgasse 7, 53111 Bonn, Germany
\endaddress
\email moree\@mpim-bonn.mpg.de \endemail
\abstract
We extend the ``character sum method" for the computation
of densities in Artin primitive root problems given by 
Lenstra and the authors [\LMS]
to the situation of radical extensions of arbitrary rank.
Our algebraic set-up identifies the key parameters of the situation at
hand, and obviates the lengthy analytic multiplicative number theory 
arguments that used to go into the computation of actual densities.
It yields a conceptual interpretation of the formulas obtained, and 
enables us to extend their range of application in a systematic way.
\endabstract
\endtopmatter

\document

\head 1. Introduction
\endhead

\noindent
\global\baselineskip=13pt
Artin's classical 1927 conjecture provides, for 
an integer $a\ne 0, \pm 1$, a value for the density of the set
of primes $q\nmid a$ for which $a$ is a primitive root modulo $q$.
For $q\nmid 2a$, the index $[\FF_q^*:\langle a\mod q\rangle]$
is divisible by $n$ if and only if $q$ splits completely 
in the splitting field $F_n$ of the polynomial $X^n-a$ over $\Que$,
so $a$ is a primitive root modulo~$q$ if and only if $q$ does
not split completely in any of the fields $F_p$, with $p<q$ prime.
One may therefore expect, with Artin, that a fraction
$$
\prod_{p\ \text{prime}} \Bigl(1-{1 \over [F_p:\Que]}\Bigr)
\tag{1.1}
$$
of all primes has $a$ as a 
primitive root.
Two problems arise when proving this.

The first problem is of an algebraic nature, and was overlooked 
for more than 30 years [\ST].
This problem, which is at the heart of this paper and
makes that 1.1 is not in general correct,
is the possibility of a {\it dependency\/} 
between the splitting conditions in the various fields $F_p$.
It arises if the Galois group of the compositum of the 
fields $F_p$ over $\Que$ is a strict subgroup of the product
group $\prod_p \Gal(F_p/\Que)$, in which it is naturally contained.
It leads to the question of how the probability for a
prime $q$ to have a certain splitting behavior in the compositum
is related to the probabilities of having prescribed splitting
behavior in each of the fields $F_p$.

The second problem is of an analytic nature.
It is caused by the fact that the well-known theorem that a `fraction'
$[F_p:\Que]^{-1}$ of all primes splits completely in $F_p$
is an asymptotic statement, and combining
these statements for all primes $p$ into a single asymptotic
statement is far from trivial, and requires good control of the error
terms in the statements involved.
Such error terms are only available under assumption of the
generalized Riemann hypothesis (GRH) for the fields $F_n$.

Hooley [\HO] showed that, when the 
algebraic problem is taken into account, the conjectural densities
provided by Artin's argument are, under GRH,
indeed the correct densities.
His argument was extended by Cooke and Weinberger [\CW] to deal, under GRH,
with the case of Artin-type densities over arbitrary number fields.

In the case of Artin's original conjecture, Hooley's argument
proves that the density of the set of primes that do not split completely
in any of the fields $F_n$ with $n>1$ equals
the inclusion-exclusion value
$$
\sum_{n=1}^\infty  {\mu(n)\over [F_n:\Que]},
\tag{1.2}
$$
with $\mu$ the M\"obius function.
For the correct evaluation of this sum, one has to control the
algebraic problem that the degree $[F_n:\Que]$ for squarefree
values of~$n$ may not equal the product
of the degrees $[F_p:\Que]$ with $p|n$.
If it does so for all squarefree $n$, then 1.2 equals the Euler product~1.1.
In general, one needs to multiply~1.1 by a rational correction factor
that takes the entanglement between the fields $F_p$ into account~[\ST].

Well-known variants of Artin's conjecture ask for the density
of primes $q$ satisfying more complicated conditions than just
having $a$ as a primitive root.
One may for instance require that the 
index of $\langle a\mod q\rangle$ in
$\FF_q^*$ be equal to or divisible by some integer, look at
$\langle a\mod q^k\rangle\subset (\Zee/q^k\Zee)^*$ for $k>1$,
or impose an additional congruence condition on $q$.
All these variants can be phrased in terms of the splitting 
behavior of $q$ in the union
$F_\infty=\bigcup_{n\ge1} F_n$ of all $F_n$
inside some algebraic closure $\overline\Que$ of~$\Que$.
The field $F_\infty$ is
obtained by adjoining to $\Que$ the group of radicals
$$
R_\infty=
\{x\in\overline\Que^*: x^k\in \langle a\rangle
\quad\text{for some $k>0$}\}.
$$
The associated Artin-densities can now be given by infinite
sums generalizing 1.2 as in [\LEa, formula 2.15],
but their explicit evaluation as a rational multiple of
the naive Euler product replacing 1.1 tends to lead to 
nasty calculations, as the field degrees encountered do not always
admit an easy description.

The problems arise from the fact that
the Galois group of $G=\Gal(\Que(R_\infty)/\Que)$
is usually not the product over $p$ of the Galois groups
$\Gal(\Que(R_{p^\infty})/\Que)$ coming from the adjunction of all
$p$-power roots.
It is however shown in [\LMS, Section 2] that
$G$ is a closed subgroup of index two of the profinite group
$$
A=\Aut_{R_\infty\cap\Que^*} (R_\infty)
$$
of group automorphisms of $R_\infty$ fixing the rational
numbers in $R_\infty$, and that this larger group $A$
does admit a natural splitting as a product $\prod_p A_p$
of automorphism groups of $p$-power radicals.

As a profinite group, $A$ naturally comes with a Haar measure~$\nu$.
The {\it character sum method\/} [\LMS, Theorem 3.3] shows that,
starting from any reasonable subset 
$$
S=\prod_p S_p\subset \prod_p A_p =A
$$
that characterizes `good splitting' at each of the $p$-components,
we can {\it always\/} decompose the associated Artin-density
$\delta(S)=\nu(G\cap S)/\nu(G)$ in a natural way
as a product
$$
\delta(S)={\nu(G\cap S)\over \nu(G)}=
E\cdot{\nu(S)\over \nu(A)}
$$
of the {\it naive density\/} $\nu(S)/\nu(A)$ and
a rational {\it entanglement correction factor\/}
$$
E=1+\prod_p E_p.
$$
Here $E_p$ denotes the average value on $S_p$ of the $p$-component
$\chi_p$ of the character
$$
\chi=\prod_p\chi_p: A=\prod_p A_p\tto\{\pm1\}
$$
that has $G$ as its kernel.
As we have $\chi_p=1$ for almost all~$p$, almost all $E_p$ equal 1,
and the correction factor $1+\prod_p E_p$ only involves a
finite number of {\it critical\/} primes~$p$ at which the
entanglement of $p$-power radicals takes place.

When applicable, the method of [\LMS] leads to smooth 
computations of the densities involved.
It is however too restrictive to deal with generalizations
of Artin's conjecture that refer to
properties of more than a single integer or rational number 
modulo the prime numbers $q$ under consideration.
Already for two rational numbers $a_1, a_2\in\Que^*$, one may 
wonder for which fraction of the primes $q$
\smallskip
\itemitem{1.}
the subgroup generated by $a_1$ and $a_2$ modulo $q$
equals $\FF_q^*$;
\itemitem{2.}
each of $a_1$ and $a_2$ is a primitive root modulo $q$;
\itemitem{3.}
$a_1$ and $a_2$ generate the same subgroup modulo $q$;
\itemitem{4.}
$a_1$ is in the subgroup generated by $a_2$ modulo $q$.
\smallskip\noindent
The first three problems immediately generalize to $k\ge 2$
rational numbers, and it is clear that many variations
of these higher rank primitive root problems exist.

In this paper, we present the extension of our rank-1 result in [\LMS]
to arbitrary ranks, and show that it leads to a simple, unified
approach to compute higher rank primitive root densities.
Section 2 provides a description of the Galois group $G$ as a subgroup
of $A$, and Section 3 achieves the decomposition of the Artin
density $\delta(S)$ given above.
With minimal effort, we then recover the Artin densities for higher
rank subgroups due to Cangelmi-Pappalardi [\CP] in Section 4,
and the multiple primitive root densities due to Matthews [\MA]
in Section 5.
Our approach explains the structure of the formulas
that are found in these papers after cumbersome manipulations
of double and triple sums arising from analogues of 1.2.
In particular, we see which special cases lead to
nice `multiplicative formulas', and show how a direct application
of the method also yields the extension of Matthews result
that was recently found by Schinzel [\SC].
A final section adresses the vanishing problem for primitive
root densities, which is more complicated in this higher rank case
than it is in the rank-1 case.


\head 2. Radical extensions of the rational number field
\endhead

\noindent
Let $\Gamma\subset\Que^*$ be a finitely generated subgroup.
We want to explicitly describe the Galois group $G$ of
the radical extension $\Que\subset\Que(\Gamma_\infty)$ obtained by 
adjoining to $\Que$ the group
$$
\Gamma_\infty=\{x\in\overline\Que^* \subset \CC^*: x^n\in \Gamma
    \text{ for some $n\in\Zee_{\ge1}$}\}
$$
of {\it complex\/} roots of arbitrary order of the elements in $\Gamma$.
The group $\Gamma_0= \Gamma_\infty\cap\Que^*$ is a finitely
generated subgroup of $\Que^*$ that contains $\Gamma$ as a subgroup
of finite index. 
We can choose a $\Zee$-basis $\{b_i\}_{i=1}^r$ of positive 
rational numbers for the subgroup
$\Gamma_0^+=\Gamma_0\cap \Que_{>0}\subset\Que^*$ of positive
elements in $\Gamma_0$ and write
$$
\Gamma_0= \Gamma_0^+\times \langle-1\rangle =
          \langle b_1\rangle \times \langle b_2\rangle \times \ldots
          \langle b_r\rangle \times \langle-1\rangle,
$$
with $r\in\Zee_{\ge0}$ the {\it rank\/} of $\Gamma$.
The subgroup of positive real numbers in $\CC^*$
naturally forms a $\Que$-vector space,
so if we denote by $\Gamma_\infty^+$ the $\Que$-vector space
generated by~$\Gamma_0^+$, we have
$$
\Gamma_\infty = \Gamma_\infty^+\times \mu_\infty 
              = b_1^\Que \times b_2^\Que \times \ldots
                b_r^\Que \times \mu_\infty,
\tag{2.1}
$$
with $\mu_\infty$ the group of roots of unity in $\CC^*$.
As in [\LMS], we analyze the Galois group 
$G=\Gal(\Que(\Gamma_\infty)/\Que)$ in terms of the 
{\it Galois representation\/}
$$
G=\Gal(\Que(\Gamma_\infty)/\Que)\longrightarrow
A=
\Aut_{\Gamma_\infty\cap\Que^*} (\Gamma_\infty)=
\Aut_{\Gamma_0} (\Gamma_\infty)
\tag{2.2}
$$
describing the action of $G$ by group automorphisms on the group
of radicals $\Gamma_\infty$.
An automorphism $\sigma\in\Aut(\Gamma_\infty)$ that leaves $\mu_\infty$
and $\Gamma_0$ invariant 
is determined by the sequences of roots of unity
$$
\left\{{\sigma(b^{1/n})\over b^{1/n}}\right\}_{n=1}^\infty
\in\prod_{n=1}^\infty \mu_n
$$ 
by which it multiplies the $n$-th roots of the elements $b\in \Gamma_0^+$
for $n\in\Zee_{\ge1}$.
Here~$\mu_n$ denotes the group of $n$-th roots of unity in 
$\CC^*$.
We can naturally view such sequences as elements of
the Tate module $\muhat=\lim\limits_{\leftarrow n} \mu_n$
of the multiplicative group, and the automorphism $\sigma$
as a $\muhat$-valued homomorphism on $\Gamma_0^+$.
In this way, we can describe~$A$ by the split
exact sequence of topological groups forming the lower row of the
diagram below:
$$

\def\mapdown#1{\Big\downarrow\rlap{$\vcenter{\hbox{$\scriptstyle#1$}}$}}
\matrix
1&\mapright{}&\Gal(\Que(\Gamma_\infty)/\Que(\mu_\infty))&\mapright{}&G
             &\mapright{}&\Gal(\Que(\mu_\infty)/\Que)&\mapright{}&1\cr
        &&\mapdown{}&&\mapdown{(2.2)}&&\mapdown{\wr}&&\cr
        1&\mapright{}&\Hom(\Gamma_0^+,\muhat)&\mapright{}&A
             &\mapright{}&\Aut(\mu_\infty)&\mapright{}&1\rlap{.}\cr
\endmatrix
$$
In order to describe $G$ as a subgroup of $A$, we consider the upper row
of the diagram, which is exact by Galois theory.
As 2.2 induces the familiar isomorphism 
$$
\Gal(\Que(\mu_\infty)/\Que)\isar\Aut(\mu_\infty)
              \iso\widehat\Zee^*=\lim_{\leftarrow n} (\Zee/n\Zee)^*
$$
occurring as the right vertical arrow in the diagram, we can
describe $G\subset A$ by identifying the image of 
$\Gal(\Que(\Gamma_\infty)/\Que(\mu_\infty))$
in $\Hom(\Gamma_0^+,\muhat)$ under the restriction of 2.2,
i.e., under the natural left vertical arrow that makes our
diagram commute.
By Kummer theory, we have an isomorphism
$$
\Gal(\Que(\Gamma_\infty)/\Que(\mu_\infty)) \isar
\Hom(\Gamma_\infty^+/(\Gamma_\infty^+\cap \Que(\mu_\infty)), \mu_\infty).
$$
As in [\LMS, Lemma 2.3], all roots of elements in $\Gamma_0^+$ that
are in the maximal abelian extension $\Que(\mu_\infty)$ of $\Que$
are {\it square\/} roots, so 
$\Gamma_\infty^+\cap \Que(\mu_\infty)=(\Gamma_0^+)^{1/2}$
is the group of positive square roots of elements in $\Gamma_0^+$.
It follows that under 2.2,
the subgroup
$\Gal(\Que(\Gamma_\infty)/\Que(\mu_\infty))\subset G$
corresponds to the subgroup 
$$
\Hom(\Gamma_\infty^+/(\Gamma_0^+)^{1/2},\mu_\infty) =
\Hom((\Gamma_0^{1/2})^+, \muhat) =
2\cdot \Hom(\Gamma_0^+, \muhat)\subset A.
$$
As $\Gamma_0^+$ is a free abelian group of rank $r$,
the group $\Hom(\Gamma_0^+, \muhat)\iso \bigoplus_{i=1}^r \muhat$ is a
free $\widehat\Zee$-module of rank $r$. 
It contains
$2\cdot \Hom(\Gamma_0^+, \muhat)$ as a subgroup of index~$2^r$.
It follows that $G\subset A$ is a subgroup of index $2^r$ as well.

In field theoretic terms, the index $2^r$ of $G$ in $A$
reflects the fact that $\Que(\Gamma_\infty)$ is the compositum of 
$\Que(\Gamma_\infty^+)$ and $\Que(\mu_\infty)$, and that
these subfields intersect in the 
`multiquadratic' extension
$\Que((\Gamma_0^{1/2})^+)=\Que(\{\sqrt {b_i}\}_{i=1}^r)$
of degree~$2^r$.
Thus, if we write an automorphism
$$
\alpha\in A= \Hom(\Gamma_0^+, \muhat)\rtimes \Aut(\mu_\infty)
$$
as $\alpha=(\phi, \sigma)$, and denote the composition of
$\phi: \Gamma_0^+\to \muhat$ with the natural map $\muhat\to\mu_2$
by~$\phi_2$, we have
$$
\alpha=(\phi, \sigma)\in G \Longleftrightarrow
\phi_2(b)= (b^{1/2})^{\sigma-1}\in\mu_2
\qquad
\text{for all $b\in \Gamma_0^+$}.
$$
Here $\sigma\in \Aut(\mu_\infty)$ acts on
$b^{1/2}=\sqrt b\in \Que(\mu_\infty)$ under the identification
$\Aut(\mu_\infty)=\Gal(\Que(\mu_\infty)/\Que)$.
We find that $G\subset A$ is the subgroup of $A$ that is
`cut out' by quadratic characters $\chi_{[b]}: A\to\mu_2$
coming from elements $b\in \Gamma_0^+$, with
$$
\chi_{[b]}: \alpha=(\phi, \sigma)\longmapsto
\phi_2(b)\cdot (b^{1/2})^{\sigma-1}.
\tag{2.3}
$$
In terms of the generators $b_1, b_2, \ldots, b_r$
of $\Gamma_0^+=\Gamma_\infty\cap \Que_{>0}$, the injective
Galois representation 2.2 fits in an exact sequence
$$
1\to G=\Gal(\Que(\Gamma_\infty)/\Que)
 \mapright{(2.2)} A
 \mapright{\bigoplus_{i=1}^r \chi_i} \mu_2^r\to 1,
\tag{2.4}
$$
where $\chi_i=\chi_{[b_i]}$ denotes
the quadratic character corresponding to $b=b_i$ as in 2.3.
In words, the exactness of 2.4 means that a group automorphism in
$A=\Aut_{\Gamma_0} (\Gamma_\infty)$
is a field automorphism in $G=\Gal(\Que(\Gamma_\infty)/\Que)$
if and only if its action on the group
$$
(\Gamma_0^+)^{1/2}=
\{x\in \overline\Que^* : x^2\in \Gamma_0^+\}=
\Gamma_\infty^+\cap\Que(\mu_\infty)
$$
induced by the inclusion map $(\Gamma_0^+)^{1/2}\subset \Gamma_\infty$
coincides with the action via the cyclotomic restriction map
$A\to \Aut(\mu_\infty)=\Gal(\Que(\mu_\infty)/\Que)\to
\Gal(\Que((\Gamma_0^+)^{1/2}/\Que)$.
We summarize the discussion in the following theorem.
\proclaim{2.5. Theorem}
Let $\Gamma\subset \Que^*$ be of rank $r\ge0$, and define
$\Gamma_0^+=\Gamma_\infty\cap \Que_{>0}$ as above.
Then $G=\Gal(\Que(\Gamma_\infty)/\Que)$ is a normal subgroup
of index $2^r$ of $A=\Aut_{\Gamma_\infty\cap\Que^*}(\Gamma_\infty)$
under the embedding $2.2$, and we have a perfect pairing
$$
\eqalign{
A/G\times \Gamma_0^+/{\Gamma_0^+}^2 &\longrightarrow \mu_2 \cr
(\alpha, b)\qquad                   &\longmapsto \chi_{[b]}(\alpha)\cr
}
$$
of elementary abelian $2$-groups defined by $2.3$.
\qed
\endproclaim
\noindent
\medskip\noindent
{\bf 2.6. Remark.}
Note that, even though we chose the generating
elements $b_i\in\Que^*$ to be positive for the sake of an easy
splitting in 2.1, the characters $\chi_{[b]}$ in 2.3
are unchanged if we replace $b$ by $-b$.
This is because for $\alpha=(\phi, \sigma)\in A$, definition~2.3 gives us
$$
\eqalign{
\chi_{[-b]}(\alpha) &= \phi_2(-b)\cdot {\sigma(\sqrt{-b})\over \sqrt{-b}}
 = \phi_2(-1)\cdot {\sigma(\sqrt{-1})\over \sqrt{-1}}\cdot \chi_{[b]}(\alpha)\cr
 &= \Bigl({\alpha(\sqrt{-1})\over\sqrt{-1}}\Bigr)^2\cdot\chi_{[b]}(\alpha)
           =\chi_{[b]}(\alpha).\cr
}
$$

\head 3. Higher rank entanglement correction
\endhead

\noindent
The group $\mu_\infty$ of roots of unity is generated
by its subgroups $\mu_{p^\infty}$ of roots of unity of
prime power order, for $p$ a prime.
In the same way,
the radical group $\Gamma_\infty$, which consists of all roots of
arbitrary order of elements in~$\Gamma$, is generated by its subgroups
$$
\Gamma_{p^\infty}= \{x\in\CC^*: x^{p^n}\in \Gamma
    \text{ for some $n\in\Zee_{\ge1}$}\}
$$
of prime power radicals. 
This gives rise to a natural
isomorphism $A\isar \prod_p A_p$, with
$A_p=\Aut_{\Gamma_{p^\infty}\cap\Que^*}(\Gamma_{p^\infty})$.
In terms of the description of $A$ provided in the previous section,
this easily follows from the decompositions
$\muhat=\prod_p \muhat_p$, with
$\muhat_p=\lim\limits_{\leftarrow n} \mu_{p^n}$, and
$\Aut(\mu_\infty)=\prod_p \Aut(\mu_{p^\infty})$.
These yield
$$
A\iso \Hom(\Gamma_0^+, \muhat)\rtimes \Aut(\mu_\infty)
= \prod_p
\left[
\Hom(\Gamma_0^+, \muhat_p)\rtimes \Aut(\mu_{p^\infty})
\right].
$$
Each of the characters $\chi_{[b]}: A \to \{\pm 1\}$ in 2.3 is continuous,
and can uniquely be written as a finite product
$\chi_{[b]}=\prod_p \chi_{[b], p}$
of {\it $p$-primary\/} quadratic characters
$$
\chi_{[b], p}: A\to A_p \to \{\pm 1\}
$$
that factor via a $p$-component $A_p$ of $A$ for some prime $p$.
This is because $\chi_{[b]}$ is defined as the product of two
quadratic characters that each have this property.

The first of these characters maps $\alpha=(\phi, \sigma)\in A$
to $\phi_2(b)$. 
It describes the action of $\alpha\in A$ on
$\sqrt b\in \Gamma_{2^\infty}\subset\Gamma_\infty$ and factors
via $A_2$.

The second character maps $\alpha=(\phi, \sigma)$ to 
$(b^{1/2})^{\sigma-1}=\sigma(\sqrt b)/\sqrt b$.
It factors via the cyclotomic component 
$\Aut(\mu_\infty)=\Gal(\Que(\mu_\infty)/\Que)\iso\widehat\Zee^*$ of $A$,
and is the lift 
$\chi_K: A \to \widehat\Zee^*\to\{\pm1\}$
of the Dirichlet character on $\widehat\Zee^*$ corresponding to
the quadratic field $K=\Que(\sqrt b)$ of discriminant $d(b)$.
We can decompose $\chi_K$ 
as
$$
\chi_K=\prod_{p|d(b)} \chi_{K,p},
\tag{3.1}
$$
with $\chi_{K,p}: A\to A_p\to\Zee_p^*\to\{\pm1\}$ the lift of
a quadratic Dirichlet character of $p$-power conductor dividing $d(b)$.
If $p|d(b)$ is odd, then $\chi_{K,p}$ is 
the lift of the Legendre symbol at $p$.
If $d(b)$ is even, there are three possibilities for $\chi_{K,2}$.
It is the character corresponding to $\Que(i)$ for $d(b)\congr 4\mod 8$,
and to one of the fields $\Que(\sqrt{\pm2})$ for $8|d(b)$.

The profinite groups $A_p$ come with a Haar measure $\nu_p$,
and if we normalize these to have $\nu_p(A_p)=1$ for all $p$,
the product measure $\nu=\prod_p \nu_p$ is a normalized Haar measure
on $A$.

The `correction factors' occurring in the densities associated to
primitive root problems find their origin in the fact that the Galois
group $G\subset A$ does not in general decompose as a product
$G=\prod_p G_p$, with
$G_p=\mathop{\text{im}}[G\to A_p]=\Gal(\Que(\Gamma_{p^\infty})/\Que)$.
The problem is to determine, for a measurable subset
$$
S=\prod_p S_p \subset \prod_p A_p = A
$$
that is given as a product of measurable subsets $S_p\subset A_p$,
the `fraction' or density $\delta(S)=\nu(G\cap S)/\nu(G)$ of elements of
$G$ that lie in $S$.
It turns out that this density can be written as the product of a
{\it naive density\/} $\nu(S)/\nu(A)$ that disregards the difference
between $G$ and $A$ and a well-structured
{\it entanglement correction factor\/}~$E$.

\proclaim{3.2. Theorem}
Let $G\subset A$ be the injection from $2.2$, and
$\nu=\prod_p \nu_p$ the Haar measure on $A=\prod_p A_p$.
Take $S=\prod_p S_p\subset A$ to be a product of
$\nu_p$-measurable subsets $S_p\subset A_p$ with $\nu_p(S_p)>0$.
Then we have 
$$
\delta(S)=
{\nu(G\cap S)\over \nu(G)} = E \cdot {\nu(S)\over \nu(A)},
$$
for an entanglement correction factor $E$ given by
$$
E=\sum_{\chi\in X} E_\chi = \sum_{\chi\in X} \prod_p E_{\chi, p}.
$$
Here $X=\Hom(A/G, \mu_2)$ denotes the dual group of $A/G$, and
the local correction factor
$$
E_{\chi, p}={1\over\nu_p(S_p)}\int_{S_p} \chi_p d\nu_p
$$
of $\chi$ at $p$ is
the average value on~$S_p$ of the $p$-primary component $\chi_p$ of
$\chi=\prod_p \chi_p$.
\endproclaim
\noindent
{\bf Proof.}
As $A/G$ is an elementary abelian 2-group of order $[A:G]=2^r$ by 2.5,
we can write the characteristic function of $G$ in $A$ as
$
{\tenbf 1}_G= 2^{-r} \sum_{\chi\in X} \chi
$.

We can compute $\nu(G\cap S)$ by integrating ${\tenbf 1}_G$
over $S\subset A$ with respect to $\nu$.
Assume $\nu(S)=\prod_p\nu_p(S_p)>0$, as the theorem trivially holds
in the case $\nu(S)=0$.
Using the equality $\nu(G)=2^{-r}\cdot\nu(A)$, we easily obtain
$$
{\nu(G\cap S)\over \nu(G)} =
{2^r \over \nu(A)}\int_S {\tenbf 1}_G\, d\nu =
{\nu(S)\over \nu(A)} \cdot
\sum_{\chi\in X} \left( {1\over\nu(S)} \int_S \chi d\nu \right) .
$$
Now $\nu(S)$ equals $\prod_p \nu_p(S_p)$,
and the integral of $\chi=\prod_p \chi_p$ over $S=\prod_p S_p$ is the product
of the values $\int_{S_p}\chi_p d\nu_p$ for all $p$.
The theorem follows.
\qed
\medskip\noindent
We deduce from Theorem 3.2 that the density $\delta(S)$ can vanish
for two reasons.
The `obvious' reason for vanishing is that the naive density
equals zero, i.e., that the set $S=\prod_p S_p$ we are looking at is 
a set of measure zero.
In Artin-like problems, this trivial reason mostly occurs in cases
(excluded in the theorem) where the set $S_p$ is empty
for some prime $p$, and we are trying to impose a kind of splitting 
behavior that is `impossible at $p$'.
For instance, a square $a$ will not be a primitive root modulo any odd prime
$q$ for the simple reason that every such $q$ splits completely in
the field $F_2=\Que(\sqrt a)=\Que$ at $p=2$.
The other, much more subtle reason is that even though $S$ has positive
measure, the splitting behavior at $p$ encoded in the sets $S_p$
is incompatible with the entanglement between the fields 
$\Que(\Gamma_{p^\infty})$ at {\it different\/} $p$.
We will discuss this further in Section 6.

\head 4. Artin's conjecture for higher rank subgroups
\endhead

\noindent
Let $\Gamma\subset \Que^*$ be a finitely generated subgroup 
of {\it positive\/} rank $r>0$.
Then for all but finitely many primes $q$, the group $\Gamma$ consists
of $q$-units, and one may ask for the density of the set of primes $q$
for which the reduction map $\Gamma \to \FF_q^*$ is surjective.
If $\Gamma=\langle r\rangle$ is cyclic of rank 1,
we are in the case of Artin's classical 1927 conjecture, and this
generalization is the most obvious higher rank analogue.
It may be analyzed in a similar way, as 
we are now after the density of the set of primes $q$ that do not
split completely in any of the fields $M_p=\Que(\Gamma^{1/p})$
generated by
$$
\Gamma^{1/p}= \{ x\in\CC^*: x^p\in\Gamma \},
$$
for $p<q$ prime.
In this case, the associated set $S=\prod_p S_p\subset A$
of `good' Frobenius elements is obtained by taking
$$
S_p=A_p\setminus \ker\varphi_p
\tag{4.1}
$$
equal to the complement of the
kernel of the natural restriction map
$$
\varphi_p: A_p \tto A(p)=\Aut_{\Gamma^{1/p}\cap\Que^*}(\Gamma^{1/p}).
\tag{4.2}
$$
The group $A(p)=\Aut_{\Gamma^{1/p}\cap\Que^*}(\Gamma^{1/p})$ is an
extension of $\Aut(\mu_p)\iso\FF_p^*$ by the dual group
$\Hom(\overline \Gamma_p,\mu_p)$ of
$$
\overline \Gamma_p=\im [\Gamma\to\Que^*/{\Que^*}^p] .
$$
As the natural map $\Que^*/{\Que^*}^p\to\Que(\mu_p)^*/{\Que(\mu_p)^*}^p$
is injective for all primes $p$, we have a natural isomorphism
$$
\Gal(M_p/\Que)=\Gal(\Que(\Gamma^{1/p})/\Que)\isar A(p)
\tag{4.3}
$$
for the Galois group of $M_p=\Que(\Gamma^{1/p})$ over $\Que$
induced by our fundamental map 2.2.
In particular, $\ker\varphi_p$ has measure 
$\nu_p(\ker\varphi_p)= (\# A(p))^{-1}=[M_p:\Que]^{-1}$
for all~$p$. 
For the measure of $S$ we find the analogue
$$
\nu(S)={\nu(S)\over \nu(A)}=
\prod_p \Bigl(1-{1\over [M_p:\Que]}\Bigr)=
\prod_p \Bigl(1-{1\over (p-1)\#\overline\Gamma_p}\Bigr)
\tag{4.4}
$$
of the naive density in 1.1.
As in the rank-1 case, the density vanishes if and only if 
${\overline \Gamma_2}$ is the trivial group, i.e.,
if and only if $\Gamma$ consists of squares in $\Que^*$.
As $\overline \Gamma_p$ has order $p^r$ for almost all primes $p$,
the density in 4.4 is a rational multiple of the
{\it rank-$r$ Artin constant\/}
$$
C_r=\prod_p \left(1-{1\over (p-1)p^r}\right).
\tag{4.5}
$$
It is a straightforward application of 3.2 to determine the
density of $G\cap S$ in $G=\Gal(\Que(\Gamma_\infty)/\Que)$.
Under GRH, this is the density of set of the primes $q$ for which
$\Gamma \to \FF_q^*$ is surjective.
\proclaim{4.6. Theorem}
Let $\Gamma\subset \Que^*$ be finitely generated and of positive rank.
Then the density inside $A$ of the subset $S=\prod_p S_p\subset A$
defined by~$4.1$ is given by $4.4$.
If\/ $\Gamma$ is not contained in ${\Que^*}^2$, then $S$ is non-empty,
and its density inside the Galois group
$G=\Gal(\Que(\Gamma_\infty)/\Que)\subset A$
from $2.2$ equals
$$
{\nu(G\cap S)\over \nu(G)}=
\Biggl(
1+\sum_{b\in{\overline \Gamma_2}\setminus \{1\} \atop d(b)\congr 1\mod 4}
  \prod_{p|2\cdot d(b)} {-1\over [M_p:\Que]-1}
\Biggr)\cdot {\nu(S)\over \nu(A)}.
$$
Here we write $M_p=\Que(\Gamma^{1/p})$ and
${\overline\Gamma_p}=\im[\Gamma\to\Que^*/{\Que^*}^p]$ as above,
and denote the discriminant of $\Que(\sqrt b)$ for 
$b\in\Que^*/{\Que^*}^2$ by $d(b)$.
\endproclaim
\noindent
{\bf Proof.}
We computed the density $\nu(S)/\nu(A)$ of $S$ in $A$ in 4.4.
For the other statement, which trivially holds for $\nu(S)=0$,
we assume that $\nu(S)=\prod_p \nu_p(S_p)$ is positive.
In this case, Theorem 3.2 expresses the density $\nu(G\cap S)/\nu(G)$ of
$G\cap S$ in $G$ as the product of $\nu(S)/\nu(A)$ and
a correction factor $E=\sum_{\chi\in X}\prod_p E_{\chi, p}$.

To see that the factor $E$ has the form stated,
we use the explicit description of $X\iso \Gamma_0^+/{\Gamma_0^+}^2$
provided by Theorem 2.5, and
compute $E_\chi=\prod_p E_{\chi, p}$ for each of the characters
$\chi=\chi_{[b]}$ with $b\in \Gamma_0^+/{\Gamma_0^+}^2\subset 
\Que^*/{\Que^*}^2$.

As $S_p=A_p\setminus \ker\varphi_p$ is the set-theoretic difference
of a group and a subgroup, the local correction factors
$$
E_{\chi, p}=
{1\over\nu(S_p)} \left[\int_{A_p} \chi_p\,d\nu_p-
                       \int_{\ker\varphi_p} \chi_p\,d\nu_p\right]
$$
at the characters $\chi$ come in three different kinds.
If $\chi_p$ is trivial, we have $E_{\chi,p}=1$.
If $\chi_p$ is non-trivial on $\ker\varphi_p$, and consequently on $A_p$,
we have $E_{\chi,p}=0$ as~both integrals vanish, being integrals of
a non-trivial character over a group.
The most interesting is the remaining third case,
in which $\chi_p$ is trivial on $\ker\varphi_p$
but not on~$A_p$.
It leads to
$$
E_{\chi, p}={-\nu_p(\ker\varphi_p)\over\nu_p(S_p)}
    ={-[M_p:\Que]^{-1}\over 1-[M_p:\Que]^{-1}}
    ={-1\over [M_p:\Que]-1}.
\tag{4.7}
$$

For our correction factor $E=\sum_{\chi\in X} E_\chi$, we need to
sum over $\chi=\chi_{[b]}$ with $b\in \Gamma_0^+/{\Gamma_0^+}^2$.
For $b=\overline 1$, we have the trivial character and obtain 
a term $E_\chi=1$.
For $b\ne \overline 1$, the field $\Que(\sqrt b)$ is real quadratic,
and $\chi$ has non-trivial components at
the primes dividing $2\cdot d(b)$.
For odd primes $p|d(b)$, this component is the lift of the
Legendre symbol at~$p$,
which is trivial on $\ker\varphi_p$, and $E_{\chi, p}$ is as in 4.7.

For $p=2$ and $b\ne \overline 1$ the situation is more involved.
Here $\chi_2$ is the product of the character
$\psi_{\sqrt b}: \alpha=(\phi,\sigma)\mapsto \phi_2(b)$,
describing the action
of $\alpha\in A$ on $b^{1/2}\in\Gamma_\infty$ as in 2.3,
and the lifted Dirichlet character $\chi_{\Que(\sqrt b), 2}$ 
of 2-power conductor from 3.1.
Note that $\chi_2$ is non-trivial, and that by Remark 2.6, we have
$$
\chi_2=\psi_{\sqrt b}\cdot \chi_{\Que(\sqrt b), 2} =
      \psi_{\sqrt {-b}}\cdot \chi_{\Que(\sqrt {-b}), 2}.
\tag{4.8}
$$
We find $E_{\chi, 2}=0=E_\chi$ in case $\chi_2$ is non-trivial on
the kernel of
$$
\varphi_2: A \tto \Aut_{\Gamma^{1/2}\cap\Que^*} (\Gamma^{1/2}) =
                  \Hom(\overline\Gamma_2, \mu_2),
\tag{4.9}
$$
so a non-trivial character $\chi\in X$ only contributes to the sum
$E=\sum_\chi E_\chi$, with value 
$$
E_\chi=\prod_p E_{\chi, p} = \prod_{p|2\cdot d(b)} {-1\over [M_p:\Que]-1},
$$
in the cases where its 2-component $\chi_2$ in 4.8 factors via the map 
$\varphi_2$ in 4.9.
As the restriction of $\varphi_2$ to the subgroup 
$\Aut(\mu_\infty)=1\rtimes \Aut(\mu_\infty)\subset A$ factors
via $\Aut(\mu_4)$, this does not happen if $d(b)$ (and therefore
$d(-b)$) is divisible by 8.
If $d(b)$ is not divisible by 8, exactly one of
$\chi_{\Que(\sqrt b), 2}$ and $\chi_{\Que(\sqrt {-b}), 2}$ is trivial,
and $\chi_2$ factors via $\varphi_2$ if and only if the 
element $b'\in\{b, -b\}$ with $d(b')\congr1\mod 4$
is contained in $\overline\Gamma_2$.
Conversely, for $b\in\overline\Gamma_2$, either $b$ or $-b$
is in $\Gamma_0^+/{\Gamma_0^+}^2$.
This leads to the sum in the statement of the theorem.
\qed
\medskip\noindent
{\bf Remark.}
Theorem 4.6 was originally proved starting from 1.2, with $M_p$ in the 
place of $F_p$.
Pappalardi [\PA] first considered a special case, and dealt with the
general case together with Cangelmi~[\CP].
Their result looks slightly different, as their
``generalized Artin constant'' does not include the
factor $\nu(S_2)$ that we have in our infinite product 4.4
describing the naive density.
To see that 4.6 agrees with Theorem~1 in~[\CP], one can write
the entanglement correction factor in 4.6 as
$$
\biggl( 1 - {1 \over [M_2:\Que]} \biggr)^{-1} \cdot
\biggl(
1-{1 \over [M_2:\Que]}\sum_{b\in{\overline \Gamma_2} \atop d(b)\congr 1\mod 4}
  \prod_{p| d(b)} {-1\over [M_p:\Que]-1} \biggr).
\leqno{(4.10)}
$$
Taking the product with 4.4, the Euler factor at $p=2$
`cancels', and one is led to a definition of 4.4
without the factor $\nu(S_2)=1 - {1 / [M_2:\Que]}$.

\head 5. Multiple primitive roots
\endhead

\medskip\noindent
One may generalize Artin's conjecture in a different direction
by asking, when given a non-empty finite subset
$\{a_1, a_2, \ldots, a_n\}\subset \Que^*\setminus\{\pm1\}$,
for the density of the set of primes~$q$ for which
{\it each\/} of the $n$ elements $a_i$ is a primitive root modulo~$q$.
This can also be phrased in terms of splitting conditions on $q$ in
the field $\Que(\Gamma_\infty)$ of Section~2, with
$\Gamma=\langle a_1, a_2, \ldots, a_n\rangle\subset \Que^*$
The subgroup generated by the elements~$a_i$.
This time, the question does not only depend on $\Gamma$, but also
on the infinite cyclic subgroups 
$$
\Gamma_i=\langle a_i\rangle \subset \Gamma 
$$
generated by each of the $a_i$.
For each $\Gamma_i$ and prime $p$, we have restriction maps
$$
\varphi_{p, i}: A_p\tto
\Aut_{\Gamma_i^{1/p}\cap\Que^*} (\Gamma_i^{1/p})
$$
as in 4.2.
As we want to determine the density of the set of primes $q$
that do not split completely in any of the fields 
$M_{i, p}=\Que(\Gamma_i^{1/p})$, with $p<q$ prime and
$i\in\{1, 2, \ldots, n\}$,
we define the set $S_p$ of good Frobenius elements at $p$ by
$$
S_p = A_p\setminus K_p\qquad \hbox{with}\quad
K_p=\bigcup_{i=1}^n\ker\varphi_{p, i},
\tag{5.1}
$$
and put $S=\prod_p S_p$ in the usual way.

As the subgroups $\ker \varphi_{p, i}$ making up $K_p$ are all contained in
the subgroup $\ker[A_p\to \Aut(\mu_p)]$, which has index $p-1$ in $A_p$,
we have $\nu_p(K_p)\le 1/(p-1)$ for all~$p$.
This shows that $S_p$ has positive measure for $p>2$.
For $p=2$, we have $K_2=A_2$ and $S_2=\emptyset$ 
if and only if there are no elements $\alpha\in A_2$ with
$\alpha(\sqrt {a_i})=-\sqrt {a_i}$ for $i=1, 2, \ldots, n$.
This occurs if and only if there is a subset
$I\subset \{1, 2, \ldots, n\}$
for which
$$
\hbox{$\prod_{i\in I} a_i\in \Que^*$ is a square}
\qquad\hbox{and}\qquad\hbox{$\#I$ is odd}.
\tag{5.2}
$$
If no such $I$ exists, 
we can uniquely define a homomorphism
$$
w_2:{\overline\Gamma_2}=\im[\Gamma\to\Que^*/{\Que^*}^2] \tto \mu_2
\tag{5.3}
$$
by putting $w_2(a_i)=-1$ for all $i$.

The intersection $\bigcap_{i=1}^n\ker\varphi_{p, i}$
is the subgroup $\ker\varphi_p$ occurring in 4.1.
The union 
$
K_p=\bigcup_{i=1}^n\ker\varphi_{p, i}
$
is rarely a subgroup for $n\ge2$, but it is a finite union of a number
$k_p$ of cosets of $\ker\varphi_p$.
In view of 4.3, the integer $k_p\ge1$
is the number of elements in $\Gal(M_p/\Que)$
that have trivial restriction to at least one
of the subfields $M_{i,p}\subset M_p=\Que(\Gamma^{1/p})$.
In terms of $k_p$, the naive density that is the analogue of~4.4
becomes
$$
\nu(S)={\nu(S)\over \nu(A)}=
\prod_p \bigl(1-\nu_p(K_p)\bigr)=
\prod_p \Bigl(1-{k_p\over [M_p:\Que]}\Bigr).
\tag{5.4}
$$
The density $\nu(S)$ vanishes if and only if 5.2 holds for
some $I$.
In the `generic' case where $\Gamma$ is freely generated of 
rank $r=n\ge1$ by the $a_i$,
and ${\overline\Gamma_p}=\im[\Gamma\to\Que^*/{\Que^*}^p]$ has
order $p^r$ for all $p$,
the integer $k_p=p^r-(p-1)^r$ equals the number of elements in
$(\Zee/p\Zee)^r$ having at least one zero-coefficient, and
the density 5.4 equals
the {\it Artin constant for $r$ primitive roots\/}
$$
D_r=
\prod_p \Bigr ( 1 - {p^r-(p-1)^r \over (p-1)p^r}\Bigr)=
\prod_p \Bigr ( 1 - {1-(1-1/p)^r \over p-1}\Bigr).
\tag{5.5}
$$
For general $\Gamma$ of free rank $r\ge1$,
the density $\nu(S)$ is a rational multiple of $D_r$.

As in the previous section,
we can apply 3.2 to find the density of $G\cap S$ in
$G=\Gal(\Que(\Gamma_\infty)/\Que)$.
Under GRH, this is the density of set of the primes $q$ for which
each of $a_1, a_2, \ldots, a_n$ is a primitive root modulo $q$.
\proclaim{5.6. Theorem}
Let $\Gamma=\langle a_1, a_2, \ldots, a_n\rangle\subset \Que^*$ be
generated by $n\ge1$ elements
$a_i\in\Que^*\setminus\{\pm1\}$.
Then the density inside $A$ of the subset $S=\prod_p S_p$
defined by~$5.1$ is given by $5.4$.
If no subset $I\subset \{1, 2, \ldots, n\}$ satisfies $5.2$,
then $S$ is non-empty, and its density inside
the Galois group $G=\Gal(\Que(\Gamma_\infty)/\Que)\subset A$
from $2.2$ equals
$$
{\nu(G\cap S)\over \nu(G)}=
\Biggl(
\sum_{a\in{\overline \Gamma_2} \atop d(a)\congr 1\mod 4} w_2(a)
  \prod_{p| d(a)} {-k_p\over [M_p:\Que]-k_p}
\Biggr)\cdot {\nu(S)\over \nu(A)}.
$$
Here $w_2$ and $k_p$ are as in $5.3$ and $5.4$, and 
we use the notation $M_p=\Que(\Gamma^{1/p})$ and
$d(a)=\text{{\rm disc}}(\Que(\sqrt a))$ as before.
\endproclaim
\noindent
{\bf Proof.}
We already computed $\nu(S)$ in 5.4, and we now apply 3.2
to find the correction factor $E=\sum_{\chi\in X} \prod_p E_{\chi,p}$.
The analysis is very similar to that in the proof of 4.6, as the
only change consists in the replacement of
$\ker\varphi_p= A_p\setminus S_p$ in 4.1 by the union
$K_p=A_p\setminus S_p$ of $k_p$ cosets of $\ker\varphi_p$ in 5.1.

For characters $\chi\in X$ for which the $p$-component
is trivial on $A_p$ or non-trivial on $\ker\varphi_p$, we find
$E_{\chi, p}=1$ and $E_{\chi, p}=0$, as before.
In particular, the characters contributing to $E$ are of the form
$\chi=\chi_{[a]}$ with $a\in\overline\Gamma_2$ satisfying
$d(a)\congr 1\mod 4$, just as in Theorem 4.6.
Let $\chi$ be such a character.
For primes $p|d(a)$, which are clearly odd, the $p$-component
$\chi_p$ of $\chi$ is the lift of the Legendre symbol. As $\chi_p$
is trivial on $K_p\subset \ker[A_p\to\Aut(\mu_p)]$, we find
$$
E_{\chi, p}={-\nu_p(K_p)\over\nu_p(S_p)}
    ={-k_p[M_p:\Que]^{-1}\over 1-k_p[M_p:\Que]^{-1}}
    ={-k_p\over [M_p:\Que]-k_p}.
$$
At $p=2$, the character $\chi_2=\psi_{\sqrt a}$ has, 
by definition of $K_2$,
the constant value $w_2(a)$ on $S_2=A_2\setminus K_2$,
so we find $E_{\chi, 2}=w_2(a)$.
The result follows. \qed
\medskip\noindent
Theorem 5.6 was originally proved by Matthews [\MA], who phrases his 
result in terms of the density
$$
c(p)= {k_p\over [M_p:\Que]}=\nu_p(K_p)
$$
of the set of primes $q\congr 1\mod p$ with the 
property that at least one of the generators $a_i$ of $\Gamma$
is a $p$-th power modulo $q$.

In the case where the set ${\Cal A}=\{a_1, a_2, \ldots, a_n\}$
consists of $n$ different prime numbers congruent to $1\mod 4$,
we have $c(p)=( 1-(1-1/p)^r )/(p-1)$,
and the sum of $2^n$ terms in the correction factor
$$
E=\sum_{a\in{\overline \Gamma_2} \atop d(a)\congr 1\mod 4} w_2(a)
  \prod_{p| d(a)} {-k_p\over [M_p:\Que]-k_p}
\leqno{(5.7)}
$$
in Theorem 5.6 can be rewritten in terms of 
$D(a)= c(a)/(1-c(a))$ as a product
$$
E = \prod_{a\in {\Cal A}} \bigl( 1+ D(a)\bigr)
  = \prod_{a\in {\Cal A}} \bigl( 1+ {c(a)\over 1 - c(a)}\bigr).
\leqno{(5.8)}
$$
For ${\Cal A}=\{a_1, a_2, \ldots, a_n\}$ an arbitrary set of $n$ 
odd prime numbers, one can adapt 5.8 to `filter out' only the terms
with discriminant $d(a)\congr 1\mod 4$ required in 5.7 by putting
$$
E = {1\over 2} \Bigl[
  \prod_{a\in {\Cal A}} \bigl( 1+ D(a)\bigr) + 
  \prod_{a\in {\Cal A}} \bigl( 1+ {\textstyle{-1\legendre a}} D(a)\bigr)
\Bigr] .
$$
Similar remarks apply to Theorem 4.6, provided that one uses the formula 4.10
or pays some attention in 4.6 to the factor at $p=2$.
It yields the formulation of the special case of Theorem 4.6 found in [\PA].

Schinzel [\SC] extends the particular case of Theorem 5.6 for 
${\Cal A}=\{a_1, a_2, \ldots, a_n\}$ a set of odd primes by
additionally imposing that the primes $q$ for which each of the
primes $a_i$ is a primitive root split in $\Que(\sqrt 2)$ and in 
$m$ additional quadratic fields of prime conductor $b\in {\Cal B}$,
with $\Cal B$ a set of primes disjoint from $\Cal A$.
The $m+1$ additional quadratic conditions change the naive density
by a factor $2^{-(m+1)}$, and the proof of his main theorem
goes through extensive calculations to obtain the correction factor
$$
E=
{1\over 2} \Bigl[
  \prod_{a\in {\Cal A}} \bigl( 1+ D(a)\bigr) 
  \prod_{b\in {\Cal B}} \bigl( 1- D(b)\bigr)  +
  \prod_{a\in {\Cal A}} \bigl( 1+ {\textstyle{-1\legendre a}} D(a)\bigr)
  \prod_{a\in {\Cal A}} \bigl( 1- {\textstyle{-1\legendre b}} D(b)\bigr)
\Bigr] .
$$
For us, it is an almost trivial modification of the previous
result.
We take $\Gamma$ generated by ${\Cal A}\cup{\Cal B}\cup\{2\}$, redefine
$S_2$ in the obvious way, and note that the homomorphism 
$w_2:{\overline\Gamma_2}=\im[\Gamma\to\Que^*/{\Que^*}^2] \tto \mu_2$
from 5.3, which is now defined by $w_2(a)=-1$ for $a\in{\Cal A}$ and
$w_2(b)=1$ for $b\in{\Cal B}\cup\{2\}$, yields again the constant value
on $S_2$ of the character $\psi_{\sqrt x}$ for $x\in \overline\Gamma_2$.
\head 6. Vanishing criteria
\endhead

\noindent
As we observed at the end of Section 3, the density $\delta(S)$
can vanish for the simple reason that $S$ is a zero-set,
or for the more subtle reason that, even though $S$ itself has positive
density, the entanglement correction factor $E$ vanishes.
The occurrence of the simple reason, which amounts to the vanishing
of the naive density, is usually easily established.
The vanishing of $E$ is only uncomplicated in the rank-1 case,
where $E=1+E_\chi=0$ implies that we have
$$
E_\chi=\prod_p E_{\chi, p}=-1
$$
for the unique non-trivial character $\chi\in X$ in Theorem 3.2.
As all $E_{\chi, p}$ are average values of characters $\chi_p$ on $S_p$,
and therefore bounded in absolute value by 1, these extreme cases
are easily found [\LMS, Corollary 3.4].
For rank $r\ge2$ however, the vanishing of the correction factor
$$
E= 1+ \sum_{\chi\in X\setminus\{1\}} E_\chi
$$
with $2^r-1$ non-trivial terms is not so easily established.

In the case of the higher rank Artin conjecture in Section 4,
the naive density vanishes if and only if 
$\Gamma\subset {\Que^*}^2$ consists of squares, making
$S_2$ into the empty set.
If the naive density is positive, then so is the actual
density, as the correction factor $E$ {\it never\/} vanishes
in this case.
This is immediately obvious in the rank-1 case, when we have
$E=1+E_\chi$ and $|E_\chi|<1$ is readily checked.
In higher rank cases, this is less immediate, but we know
{\it a priori\/} that the density will increase if we enlarge
$\Gamma$ from a rank-1 to a higher rank subgroup of $\Que^*$,
so there is an easy way out that bypasses an exact analysis of $E$.

For the multiple primitive root densities in Section 5,
the vanishing of the naive density is almost as simple:
it vanishes if and only if $S_2$ is empty, and this is the case
where 5.2 applies for some subset $I\subset \{1,2,\ldots, n\}$.
For the vanishing of the density in cases where the naive density
is positive, we are dealing with a rather complicated correction factor
in Theorem 5.6, and in this case the density {\it decreases\/}
if we add elements $a_i$.

In such cases, a promising way to proceed is often to {\it not\/}
directly use the formula itself, but the structural idea giving 
rise to it.
What Theorem 2.5 expresses is that the nature of the entanglement
lies in the relation the splitting behavior at 2 bears to
the splitting behavior at the {\it finitely many\/} odd primes $p$
at which the elements of $\Gamma$ can have non-zero valuation.
In cases where the splitting condition $S_p\subset A_p$ at $p$ 
factors through a finite quotient $A_p\to \bar A_p$, in the sense
that $S_p$ is the inverse image of a subset $\bar S_p\subset \bar A_p$,
it suffices to find a single element of $G$ for which the 
$p$-component at $p=2$ and at the finitely many odd critical primes
projects into $\bar S_p$ under $A_p \to \bar A_p$.
This is essentially an explicit version of the observation
in [\LEa, Theorem 4.1] that the density can only vanish due
to obstructions `at a finite level'
that can be made precise in terms of the input data.
In the examples we have in Sections 4 and 5, the map
$\varphi_p: A_p \to A(p)\iso\Gal(M_p/\Que)$ from 4.2 and 4.3
is such a finite quotient map on $A_p$ that is used to define $S_p$.

As an illustration, let us derive the vanishing conditions
for the entanglement correction factor $E$ in Theorem 5.6
without directly considering the expression for~$E$ that is given
in the theorem.
We only use the fact that the `splitting condition' $S_p$ that is
imposed reflects an actual splitting condition in the field $M_p$:
the set $S_p$ consists of Frobenius classes of primes
that do not split completely in any of the subfields
$M_{i, p}=\Que(\Gamma_i^{1/p})$.

The assumption that $S_2$ and therefore $S$ be non-empty amounts to
saying that there exist primes $q$ that are inert in all $n$ quadratic 
fields $\Que(\sqrt {a_i})$, thus satisfying the necessary splitting condition
in $M_2=\Que(\Gamma^{1/2})$.
By Theorem 2.5,
the only possible implication this can have for the splitting behavior
of $q$ in fields $M_p=\Que(\Gamma^{1/p})$ for primes $p>2$ is that
possibly, the splitting behavior of $q$ in the quadratic subfield
$\Que(\sqrt{\pm p}) \subset \Que(\mu_p)\subset M_p$ can no
longer be freely prescribed at the critical primes $p$.
For primes $p\ge5$, this will not lead to an incompatibility of
local splitting behaviors, as there will be primes $q$ with Frobenius in
$S_2$ (or even with prescribed Frobenius in $\Gal(M_2/\Que)$)
that satisfy $q\not\congr 1\mod p$ for all critical primes $p\ge5$.
Such $q$ trivially have Frobenius in $S_p$.

For $p=3$, where we have $\Que(\sqrt{-3})=\Que(\mu_3)$,
it may however happen that $-3\in \Que^*/(\Que^*)^2$
is in $\overline\Gamma_2$, and that we have $w_2(-3)=1$.
Then {\it every\/} $q$ with Frobenius in $S_2$ is congruent to $1\mod 3$,
and it can only have Frobenius class in $S_3$ if there exists a character
$$
\chi:{\overline\Gamma_3}=\im[\Gamma\to\Que^*/{\Que^*}^3] \tto \mu_3
\leqno{(6.1)}
$$
that is non-trivial on $a_i\mod {\Que^*}^3$ for $i=1,2,\ldots n$.
This yields the following algebraic (and more precise) version of the
vanishing result for multiple primitive roots. Its derivation
in [\MA, p. 114--118 and p. 138--145] takes a dozen of pages.
\proclaim{6.2. Theorem}
Let $\Gamma=\langle a_1, a_2, \ldots, a_n\rangle\subset \Que^*$ and
$S$ be as in Theorem $5.6$, and
suppose no subset $I\subset \{1, 2, \ldots, n\}$ satisfies $5.2$.
Then $S\cap \Gal(\Que(\Gamma_\infty)/\Que)$
has zero density if and only if the following conditions are satisfied:
\smallskip
\itemitem{\rm (a)}
the kernel of the map $w_2: {\overline\Gamma_2} \to \mu_2$ in
$5.3$ contains $-3\mod {\Que^*}^2;$
\itemitem{\rm (b)}
the kernel of every character $\chi:{\overline\Gamma_3}\to\mu_3$ contains
at least one element $a_i\mod {\Que^*}^3$.
\smallskip\noindent
If none of the $a_i$ is a cube in $\Que^*$, then
condition {\rm (b)} does not hold if we have $n\le3$, or
if the $\FF_3$-rank of ${\overline\Gamma_3}$ is either $1$ or at least $n-1$.
\endproclaim
\medskip\noindent
{\bf Proof.}
The hypothesis concerning 5.2 means that $S$ itself has positive density.
As explained above, $\nu(S\cap G)=0$ in 5.6 can then only occur
in the case where condition (a) is satisfied, and the splitting 
condition $S_2$ implies that we have primes $q\congr 1\mod 3$ only.
For such $q$, the condition at 3 that $q$ needs to satisfy is
that no Frobenius element over $q$ in $\Gal(M_3/\Que(\mu_3))$ fixes
a cube root $a_i^{1/3}$,
as this is equivalent to $a_i$ being a cube modulo $q$.
The group of characters $\chi:{\overline\Gamma_3}\to\mu_3$
in 6.1 may be identified, by Kummer theory and the injectivity of
the map $\Que^*/{\Que^*}^3 \to \Que(\mu_3)^*/{\Que(\mu_3)^*}^3$,
with the Galois group $\Gal(M_3/\Que(\mu_3))$.
Condition (b) therefore amounts to saying that every element of
$\Gal(M_3/\Que(\mu_3))$ fixes a cube root $a_i^{1/3}$, and we
obtain the first half of our theorem.

We finally have to deal with the question whether there exists a
character $\chi:{\overline\Gamma_3}\to\mu_3$ that is 
non-trivial on each of the generators $\overline a_i=a_i\mod {\Que^*}^3$ of
$\overline\Gamma_3$.
For this, it is clearly necessary that each
$\overline a_i$ is not the trivial element, i.e.,
none of the $a_i$ is a rational cube.
This is also sufficient if the $\FF_3$-rank of ${\overline\Gamma_3}$ is $n$,
since then a suitable character $\chi$ can simply be defined by choosing
non-trivial values $\chi(\overline a_i)$.
If the $\FF_3$-rank equals $n-1>0$, there is a single relation
expressing one generator, say $\overline a_n$, as a product of
some other generators $\overline a_i$ with exponents $\pm1$.
We now choose the values $\chi(\overline a_i)\in\mu_3\setminus\{1\}$
for $1\le i<n$ such that we have $\chi(\overline a_n)\ne 1$.
For $\Gamma_3$ of $\FF_3$-rank 1, any isomorphism
$\chi:{\overline\Gamma_3} \isar \mu_3$ does what we want.

For $n\le3$, we are automatically in one of the three cases 
we just dealt with.
\qed
\medskip\noindent
For $n=4$ and ${\overline\Gamma_3}$ of rank 2, it is possible 
that both conditions of the theorem hold without any $a_i$
being a cube. This follows from the fact that
if $a_1$ and $a_2$ are arbitrary and $a_3$ and $a_4$
are chosen to satisfy $\overline a_3= \overline a_1 \overline a_2$
and $\overline a_4=\overline a_1\overline a_2^{-1}\in {\overline\Gamma_3}$,
condition (b) will always be satisfied.
It is easy to satisfy the first condition as well: 
$(a_1, a_2, a_3, a_4)=(5,\  -3\cdot 5,\  2^3\cdot 3\cdot 5^2,\ 3\cdot 7^3)$
yields an example.

\head Appendix: numerical values
\endhead

\noindent
The rank-$r$ Artin constant $C_r$ from 4.5 and the 
Artin constant for $r$ primitive roots $D_r$ from 5.5
are defined by Euler-products that converge slowly, especially for
small values of $r$.
However, there are techniques for rapidly obtaining accurate approximations,
see [\MO].

The table below provides 20 digit decimal approximations of 
$C_r$ and $D_r$ for $1\le r \le 7$.
\bigskip\bigskip

\line{\hfill\vbox{
\offinterlineskip
\hrule
\tabskip=0pt
\halign{&\vrule#\tabskip.9em plus2em minus.6em&\strut\hfil$#$&\vrule#&
  \hfil$#$\hfil&\vrule#&\hfil$#$\hfil&\vrule#\tabskip=0pt\cr
height3pt&\omit&&\omit&&\omit&\cr
&r&&C_r&&D_r&\cr
height3pt&\omit&&\omit&&\omit&\cr
\noalign{\hrule}
height2pt&\omit&&\omit&&\omit&\cr
\noalign{\hrule}
height3pt&\omit&&\omit&&\omit&\cr
&1&&0.37395~58136~19202~28805&&0.37395~58136~19202~28805&\cr
&2&&0.69750~13584~96365~90328&&0.14734~94000~02001~45807&\cr
&3&&0.85654~04448~53542~17442&&0.06082~16551~20305~08600&\cr
&4&&0.93126~51841~60004~33438&&0.02610~74463~14917~70808&\cr
&5&&0.96666~88685~96777~51274&&0.01156~58420~47143~35542&\cr
&6&&0.98368~26363~12342~05850&&0.00525~17580~26977~39754&\cr
&7&&0.99195~72807~75518~31567&&0.00243~02267~63032~72703&\cr
height3pt&\omit&&\omit&&\omit&\cr
}\hrule}\hfill}

\enddocument